\documentclass[a4paper,12pt]{article}
\usepackage{amsmath}
\usepackage{amssymb}
\usepackage{tabularx}
\usepackage{enumerate}
\usepackage{graphicx}
\usepackage{subfigure}
\usepackage{color}

\newtheorem{thm}{Theorem}[section]
\newtheorem{lem}[thm]{Lemma}
\newtheorem{cor}[thm]{Corollary}

\newtheorem{rmk}{Remark}[section]
\newtheorem{defi}{Definition}[section]
\newtheorem{pppp}{Proof}

\newcommand{\qed}{\hspace{1em}\mbox{\raisebox{0.65ex}{\fbox{}}}}

\numberwithin{equation}{section}
\newcommand{\ol}{\overline}

\newcommand{\be}{\begin{equation}}
\newcommand{\ee}{\end{equation}}
\newcommand\bes{\begin{eqnarray}} \newcommand\ees{\end{eqnarray}}
\newcommand{\bess}{\begin{eqnarray*}}
\newcommand{\eess}{\end{eqnarray*}}

\newcommand{\R}{\mathbb{R}}
\newcommand{\bpf}{{\bf Proof:\ \ }}
\newcommand{\epf}{\mbox{}\hfill $\Box$}

\begin{document}

\thispagestyle{empty}

\title{Spreading and vanishing in a West Nile virus model with expanding fronts\thanks{The work is supported by the NSFC of China (Grant No.11371311) and
 TAPP of Jiangsu Higher Education Institutions (No.PPZY2015B109).
}}
\date{\empty}

\author{Abdelrazig K. Tarboush$^{a,b}$, Zhigui Lin$^{a}$\thanks{Corresponding author. Email: zglin68@hotmail.com} and Mengyun Zhang$^{a}$\\
{\small $^a$School of Mathematical Science, Yangzhou University, Yangzhou 225002, China}\\
{\small $^b$Department of Mathematics, Faculty of Education,}\\
 {\small University of Khartoum, Khartoum 321, Sudan}
}
 \maketitle

\begin{quote}
\noindent
{\bf Abstract.} { 
\small In this paper, we study a simplified version of a West Nile virus model discussed by Lewis et al. \cite{LRD}, which was considered as a first approximation for the spatial spread of WNv. The basic reproduction number $R_0$ for the non-spatial epidemic model is defined and a threshold parameter $R_0 ^D$ for the corresponding problem with null Dirichlet boundary condition is introduced. We consider a free boundary problem with coupled system, which describes the diffusion of birds by a PDE and the movement of mosquitoes by a ODE. The risk index $R_0 ^F (t)$ associated with the disease in spatial setting is represented. Sufficient conditions for the WNv to eradicate or to spread are given. The asymptotic behavior of the solution to system when the spreading occurs are considered. It is shown that the initial number of infected populations, the diffusion rate of birds and the length of initial habitat exhibit important impacts on the vanishing or spreading of the virus. Numerical simulations are presented to illustrate the analytical results.
 }



\end{quote}

\section{Introduction}

Emerging and re-emerging infectious diseases are growing threats to public health, agriculture and wildlife management \cite{AM, AR, BM, JE, KH, KW, TC}.
Some threats are caused by mosquitos. There are over 2500 kinds of mosquito in the world. They can transmit disease by bacterial, viruses or parasites without being affected themselves. Diseases transmitted by mosquitoes include malaria, dengue, West Nile virus, chikungunya, filariasis, yellow fever, Japanese encephalitis, Saint Louis encephalitis, Western equine encephalitis, Venezuelan equine encephalitis, Eastern equine encephalitis and Zika fever.

 West Nile virus (WNv) is an infectious disease spreading through interacting vectors (mosquitoes) and reservoirs (birds) \cite{GKLZ}, the virus infects and causes disease in horse and other vertebrate animals; humans are usually incidental reservoirs \cite{ACT, AL}.
WNv was first isolated and identified in 1937 from the blood of a febrile ugandan woman during research on yellow fever virus \cite{ACT}. Although WNv is endemic in some temperate and tropical regions such as Africa and the Middle East, it has now spread to North America, the first epidemic case was detected in New York city in 1999  and migrating birds was blamed for this introduction \cite{ACT, LRD, N, WBL}. WNv outbroke in North America in 2012 and resulted in numerous human infections and death \cite{BGW}.

As we know, no effective vaccine for the virus is currently
available and antibiotics cannot work since a virus, not bacteria, causes West Nile disease.
Therefore no specific treatment for WNv exists other than supportive therapy for severe cases, and using of mosquito repellent
 becomes the most effective preventive measure. Mathematically, it is important to understand the transmission dynamics of WNv.
  WNv yields an opportunity to explore the ecological link between vector and reservoir species.
 Taking this ecological factor into a dynamic system allows the evaluation of several control strategies \cite{WBL}.

Mathematical compartmental models for WNv have been investigated \cite{BHM, WBL}, the studies during 1950s in Egypt and Nile delta led to great advances in understanding the ecology of WNv \cite{KS}. However, most early models have only scrutinized the non-spatial dynamical formulation of the model. In recent years, spatial diffusion has been recognized as important factor to affect the persistence and eradication of infectious disease such as measles, malaria, dengue fever and WNv.

 In 2006, Lewis et al. \cite{LRD} developed and analysed a reaction-diffusion model for the spatial spread of WNv by spatially extending the non-spatial dynamical model for cross infection between birds and mosquitoes \cite{WBL}. To utilize the cooperative nature of cross-infection dynamics and analyze the traveling wave, Lewis et al. proposed in \cite{LRD} the following simplified WNv model
\begin{eqnarray}
\left\{
\begin{array}{ll}
\frac{\partial I_b}{\partial t}=D_1\Delta I_b +\alpha_{b}\beta_{b}\frac{(N_b-I_b)}{N_b} I_m -\gamma_{b}I_b ,&(x,t)\in \Omega \times (0,+\infty), \\
\frac{\partial I_m }{\partial t}=D_2\Delta I_m +\alpha_{m}\beta_b \frac{(A_m-I_m)}{N_b} I_b -d_{m} I_m ,&(x,t)\in \Omega \times (0,+\infty),\\
I_b(x,0)=I_{b,0}(x),\ I_m(x,0)=I_{m,0}(x),&x\in \overline \Omega,
\end{array} \right.
\label{Aa1}
\end{eqnarray}
where the constants $N_b$ and $A_m$ denote the total population of birds and adult mosquitos, respectively; $I_b(x,t)$ and $I_m(x,t)$ represent the populations of infected birds and mosquitos at the location $x$ in the habitat $\Omega\subset \R^n$ and at time $t\geq 0$, respectively.
The parameters in the above system are defined as follows:

$\bullet$ $\alpha_{m}$ , $\alpha_{b}$ : WNv transmission probability per bite to mosquitoes and birds, respectively;

$\bullet$ $\beta_{b}$ : biting rate of mosquitoes on birds;

$\bullet$ $d_m $ : death rate of adult mosquitos;

$\bullet$ $\gamma_b $ : bird recovery rate from WNv.

The positive constants $D_1$ and $D_2$ are diffusion coefficients for birds and mosquitoes, respectively.
Since mosquitoes do not move quickly as birds, we natually assume that $D_2\ll D_1$.

If no diffusion (i.e. $D_1=D_2=0$), then (\ref{Aa1}) becomes the spatially-independent model,
\begin{eqnarray}
\left\{
\begin{array}{ll}
\frac{d I_b(t)}{d t}=-\gamma_{b}I_b(t)+\alpha_{b}\beta_{b} \frac{(N_b-I_b (t))}{N_b}I_m(t),&t>0, \\
\frac{d I_m(t)}{d t}=-d_{m}I_m(t)+\alpha_{m}\beta_{b} \frac{(A_m-I_m (t))}{N_b} I_b (t),&t>0.
\end{array} \right.
\label{aode1}
\end{eqnarray}
It was shown in \cite{LRD} that if $0 < R_0 (:=\frac{\alpha_{m}\alpha_{b}\beta_b ^2 A_m}{d_m\gamma_b N_b})< 1$, then the virus always vanishes, while for $R_0 > 1$, a nontrivial epidemic level appears, which is globally asymptotically stable in the positive quadrant.
For system (\ref{Aa1}), if we assume that the mosquitoes population do not
diffuse ($D_2=0$), we can introduce a threshold parameter $R_0^D(:=\frac {\alpha_{m}\alpha_{b}\beta_{b}^2 A_{m}}{N_{b} d_m (\gamma_{b}+D_1\lambda_1)})$ such that for
$0 < R_0^D < 1$, the epidemic eventually tends to extinction, while for $R_0^D> 1$, a spatially inhomogeneous stationary endemic state appears
 and is globally asymptotically stable, where $\lambda_1$ is
the first eigenvalue of the boundary value problem $-\Delta \phi=\lambda\phi$ in $\Omega$ with null Dirichlet boundary condition on $\partial \Omega$.

 It is well known that the solution to (\ref{Aa1}) with null Neumann boundary condition or with null Dichlet boundary condition is positive for any positive time, this means that the environment considered is always infected, which does not match the fact that the disease appears in a small habitat and spreads gradually to a large environment.
To describe such a gradual spreading process and the changing of the infected environment, the free boundary problem has been recently introduced in some epidemic models \cite{ABL, GKLZ, LKL} and has also successfully used in other applied areas. For example,
 the melting of ice in contact with water \cite{R}, tumor growth \cite{Tao}, wound healing \cite{CF}, information diffusion in online social networks \cite{LLW} and the spreading of invasive species \cite{DG, DL2, DB, GW,LLZ, W, WC}.

 A special case of the well-known Stefan condition was derived in \cite{LLZ} by assuming that
  the amount of invasive species moving across the boundary decides the length of the expanding interval.
  Such a free boundary condition has been successfully used in \cite{DL} to describe the spreading front of invasive species by considering the logistic problem
\begin{eqnarray}
\left\{
\begin{array}{lll}
u_{t}-d u_{xx}=u(a-b u),\; &0<x<h(t),\; t>0,   \\
u_x(0,t)=u(h(t),t)=0,& t>0,\\
h'(t)=-\mu u_{x}(h(t),t),\;  & t>0, \\
u(x,0)=u_{0}(x), & 0\leq x\leq h_0,
\end{array} \right.
\label{a91}
\end{eqnarray}
here $x = h(t)$ is the free boundary to be determined, the unknown $u(x,t)$ stands for the population density of an invasive species.

 In \cite{DL}, the spreading-vanishing dichotomy was presented. The authors showed that as time approaches to infinity, the population either spreads to all new environment and successfully establishes itself, or vanishes in the long run.
 Inspired by the above research, we are attempting to consider the following simplified WNv model with the free boundary
\begin{eqnarray}
\left\{
\begin{array}{ll}
\frac{\partial I_{b}}{\partial t}=D_1\frac{\partial^2I_b}{\partial x^2}-\gamma_{b}I_{b}+\alpha_{b}\beta_{b}\frac{(N_b -I_b )}{N_b} I_m ,&g(t)<x<h(t),\, t>0, \\
\frac{\partial I_m}{\partial t}=-d_{m}I_{m}+\alpha_{m}\beta_{b}\frac{(A_m -I_m)}{N_b} I_{b},&g(t)<x<h(t),\, t>0,\\
I_b(x,t)=I_m(x,t)=0,\, & x=g(t)\, \textrm{or}\, x=h(t),\, t\geq 0,\\
g(0)=-h_0,\; g'(t)=-\mu \frac{\partial I_b}{\partial x}(g(t), t), & t>0, \\
 h(0)=h_0, \; h'(t)=-\mu \frac{\partial I_b}{\partial x}(h(t), t), & t>0,\\
I_b(x,0)=I_{b,0}(x),\ I_m(x,0)=I_{m,0}(x),&-h_0\leq x\leq h_0,
\end{array} \right.
\label{a3}
\end{eqnarray}
where $x=g(t)$ and $x=h(t)$ are the left and right moving
boundaries to be determined, $h_0$ and $\mu $ are positive constants, $\mu$ represents the expanding capability of the infected birds, the initial functions
$I_{b,0}$ and $I_{m,0}$ are nonnegative and satisfy
\begin{eqnarray}
\left\{
\begin{array}{l}
I_{b,0}\in C^2([-h_0, h_0]),\, I_{b,0}(\pm h_0)=0\, \textrm{and} \, 0< I_{b,0}(x) \leq N_b ,\  x\in (-h_0, h_0), \\
I_{m,0}\in C^2([-h_0, h_0]), I_{m,0}(\pm h_0)=0\, \textrm{and} \ 0< I_{m,0}(x) \leq A_m ,\  x\in (-h_0, h_0).
\end{array} \right.
\label{Ae1}
\end{eqnarray}
In this paper we will focus on the expanding of the infected birds and the movement of the infected mosquitoes, and study the long time behaviors of free boundaries which describe the spreading fronts of WNv.

When we finish this manuscript, we found that the recent paper \cite{WC} considered a general degenerate reaction-diffusion system with free boundary.
The results are similar, but some techniques are different, for example, we present a new way to deal with the existence of the solution. Moreover, we consider the spreading or vanishing from the epidemic view and introduce the basic reproduce numbers and the risk index of the virus.

The rest of this paper is organized as follows. In section 2, the global existence and uniqueness of the solution to problem (\ref{a3}) are proved by the contraction mapping theorem, and the comparison principle is presented.
Section 3 is devoted to the sufficient conditions for the WNv to vanish, the basic reproduction numbers and the risk index are defined.
Section 4 deals with the spreading of WNv, the sharp threshold related to the expanding capability is given and the asymptotic behavior of the solution when spreading occurs is discussed. Some simulations and a brief discussion are given in section 5.

\section{Existence and uniqueness}

In this section, we first prove the following local existence and
uniqueness results of the solution to (1.4) by the contraction mapping theorem. We then use
suitable estimates to show that the solution is defined for all $t>0$.
\begin{thm} For any given $(I_{b,0}, I_{m,0})$ satisfying $(1.5)$, and any $\alpha \in (0, 1)$, there is a constant $T>0$ such that
problem $(1.4)$ admits a unique solution
$$(I_b, I_m; g, h)\in [C^{ 1+\alpha,(1+\alpha)/2}(D_{T})]^2\times [C^{1+\alpha/2}([0,T])]^2,$$
moreover,
\begin{eqnarray} \|I_b,\, I_m\|_{C^{1+\alpha,
(1+\alpha)/2
}({D}_{T})}+||g,\, h\|_{C^{1+\alpha/2}([0,T])}\leq
C,\label{b12}
\end{eqnarray}
where $D_{T}=\{(x, t)\in \R^2: x\in [g(t), h(t)], t\in [0,T]\}$, $C$
and $T$  depend only on $h_0, \alpha, \|I_{b,0}\|_{C^{2}([-h_0, h_0])}$ and $\|I_{m,0}\|_{C^{2}([-h_0, h_0])}$.
\end{thm}
\bpf
 The proof can be proved by the similar way to \cite{ABL} or \cite{WC} with some minor modifications.
First, we are going to use $g, h$ and $I_b$ to express $I_m$, since the second equation of the model (\ref{a3}) for $I_m$ is an ODE.
For any given $T>0$, take
$$\mathcal{G}_{T}=\{g\in C^1([0,T]): \,  g(0)=-h_0,\,  g'(t)\leq 0,\, 0\leq t\leq T\},$$
$$\mathcal{H}_{T}=\{h\in C^1([0,T]): \, h(0)=h_0,\,  h'(t)\geq 0,\, 0\leq t\leq T\}.$$
Define the extension mapping $E$ by
$E_t(w)(x,t)=w(x,t)$ when $x\in [g(t), h(t)]$, and $E_t(w)(x,t)=0$ otherwise.
If $g(t)\in \mathcal{G}_{T}$, $h(t)\in \mathcal{H}_{T}$ and $I_b (x,t)\in C(D_T)$, then
$I_m$ can be represented as
$$
I_m (x,t):=H(t,I_b (x,t))=e^{-w(x,t)}\Big( E_0(I_{m,0}) (x)+\int_0^t \frac{\alpha_m \beta_b A_m}{N_b} e^{w(x,\tau)}E_{\tau}(I_b) (x,\tau)d\tau \Big)
$$
for $(x,t)\in D_T$, where $$w(x,t)=d_{m}t+\int_0^t \frac{\alpha_{m}\beta_{b}}{N_{b}} E_s(I_{b})(x,s)ds.$$

To circumvent the difficulty induced by the double free boundaries, we next straighten them. As in \cite{ABL} and \cite{DL} (see also \cite{WC}), we  make the following change of
variable:
$$y=\frac{2h_0x}{h(t)-g(t)}-\frac{h_0(h(t)+g(t))}{h(t)-g(t)},\
u(y,t)=I_b (x,t).$$
Then problem (\ref{a3}) can be transformed into
\begin{eqnarray}
\left\{
\begin{array}{ll}
u_{t}=Au_{y}+Bu_{yy}-\gamma_{b}u+\alpha_{b}\beta_b \frac{(N_b -u)}{N_b} H(t,u(y,t)),\; &t>0, \ -h_0<y<h_0, \\
u=0,\quad h'(t)=-\frac{2h_0\mu}{h(t)-g(t)}\frac{\partial u}{\partial
y},\quad &t>0, \ y=h_0,\\
u=0,\quad g'(t)=-\frac{2h_0\mu}{h(t)-g(t)}\frac{\partial u}{\partial
y},\quad &t>0, \ y=-h_0,\\
h(0)=h_0, \quad g(0)=-h_0, &\\
u(y,0)=u_0(y):=I_{b,0}(y), \; &-h_0\leq y\leq
h_0,
\end{array} \right.
\label{Hb}
\end{eqnarray}
where $A=A(h, g,
y)=y\frac{h'(t)-g'(t)}{h(t)-g(t)}+h_0\frac{h'(t)+g'(t)}{h(t)-g(t)}$, and
$B=B(h, g)=\frac{4h_0^2D_1}{(h(t)-g(t))^2}$. After this transformation, the unknown boundaries $x=h(t)$ and $x=g(t)$ become the fixed lines
$y=h_0$ and $y=-h_0$, respectively.

Denote $g^* =-\mu I'_{b,0}(-h_0), h^* =-\mu I'_{b,0}(h_0)$, and $\Delta _T=[-h_0, h_0]\times [0, T]$. For $0<T\leq1$, set
 \begin{eqnarray*}
&\mathcal{D}_{1T}=\{u\in C(\Delta_T):\, u(0, y)=I_{b,0}(y),\, u(\pm h_0,t)=0,\,
||u-I_{b,0}||_{C(\Delta_T)}\leq 1\},\\
&\mathcal{D}_{2T}=\{g\in C^1([0,T]): \  g(0)=-h_0,\, g'(0)=g^*,\
g^*-1\leq g'(y,t)\leq 0\},\\
&\mathcal{D}_{3T}=\{h\in C^1([0,T]): \  h(0)=h_0,\ h'(0)=h^*,\
0\leq h'(y,t)\leq h^*+1\}.
\end{eqnarray*}
Owing to $g(0)=-h_0$ and $h(0)=h_0$, one can see that $\mathcal{D}:=\mathcal
{D}_{1T}\times\mathcal {D}_{2T}\times\mathcal {D}_{3T}$ is a complete metric space with the
metric
$$d((u_1, g_1, h_1), (u_2, g_2,
h_2))=||u_1-u_2||_{C(\Delta_T)}+||g'_1-g'_2||_{C([0, T])}+||h'_1-h'_2||_{C([0, T])}.$$

Next, take a mapping $\mathcal {F}:\mathcal{D}
 \rightarrow C(\Delta_T)\times C^1([0,T])\times C^1([0,T])$  by
$$\mathcal{F}(u, g, h)=(\ol u, \ol g, \ol h),$$
where $\ol u\in C^{(1+\alpha)/2,1+\alpha}(\Delta_T)$ is the unique solution of
the following initial boundary value problem
\begin{eqnarray}
\left\{
\begin{array}{lll}
\ol u_{t}=A \ol u_{y}+B\ol  u_{yy}-\gamma_{b}u(y,t)&\\
\qquad +\alpha_{b} \beta_{b}\frac{(N_b-u(y,t))}{N_b} H(t, u(y,t)),\; &t>0, \ -h_0<y<h_0, \\
\ol u(-h_0,t)=\ol u(h_0,t)=0,\quad  &t>0, \\
\ol u(y,0)=u_0(y):=I_{b,0}(y), \; &-h_0\leq y\leq h_0,
\end{array} \right.
\label{b1}
\end{eqnarray}
with
\begin{equation}
\ol g(t)=-h_0-\int^t_0 \frac{2h_0\mu}{h(\tau)-g(\tau)}\frac{\partial \ol u(-h_0,\tau)}{\partial
y}d\tau,
 \label{e2}
\end{equation}
\begin{equation}
\ol h(t)=h_0-\int^t_0 \frac{2h_0\mu}{h(\tau)-g(\tau)}\frac{\partial \ol u(h_0,\tau)}{\partial
y}d\tau.
 \label{e3}
\end{equation}
  The remainder of the proof is similar as that in \cite{ABL}, \cite{DL} and \cite{WC}.
 By applying the standard $L^p$ theory and the Sobolev
imbedding theorem, one can see that for $T>0$ small enough,
$\mathcal{F}$ maps $\mathcal{D}$ into itself and $\mathcal {F}$ is a contraction mapping on
$\mathcal{D}$. So by the contraction
mapping theorem, $\mathcal{F}$ admits a unique fixed point $(u; g, h)$
in $\mathcal{D}$. Moreover, using the Schauder's estimates yields that $(u(y,t); g(t),h(t))$ is a solution of the problem (\ref{Hb}),
in other words, $(I_b (x,t), I_m (x,t); g(t),h(t))$ is a unique
local classical solution of problem (\ref{a3}).
\epf

\bigskip

The global existence of the solution to $(1.4)$ is guaranteed by the following estimates.
\begin{lem} Let $(I_b , I_m; g, h)$ be a solution to problem  \eqref{a3} defined for $t\in (0,T_0]$ for some $T_0\in (0, +\infty)$.
Then we have
\[
0<I_b (x, t)\leq N_{b}\; \mbox{ for }\, g(t)<x<h(t),\; t\in (0, T_0], \]
\[
0<I_m (x, t)\leq A_{m}\; \mbox{ for }\, g(t)<x<h(t),\; t\in (0, T_0], \]
\[
 0<-g'(t),\ h'(t)\leq C_1 \; \mbox{ for } \; t\in (0,T_0], \]
 where $C_1$ is independent of $T_0$.\label{mono}
\end{lem}
\bpf It is easy to see that $I_b$ and $I_m$ are positive, since their initial values are nontrivial and nonnegative, and system (1.4) is
quasimonotone nondecreasing. It follows from the condition (\ref{Ae1}) that $I_b\leq N_b$ and $I_m\leq A_m$ directly.

Using the Hopf boundary lemma to the equation of $I_b$ yields that
 $$ \frac{\partial I_b }{\partial x } (h(t), t)<0 \ \;\; \textrm{for} \ 0<t\leq T_0.$$
 Hence $h'(t)>0$ for $t\in (0, T_0]$ by the free boundary condition in (\ref{a3}). Similarly, $g'(t)<0$ for $t\in (0, T_0]$.

It remains to show that $-g'(t),\ h'(t)\leq C_1$ for $t\in (0,T_0]$ and
some $C_1$. The proof is similar to that of Theorem 2.3 in \cite{ABL}, see also Lemma 2.2 in \cite{WC}, we therefore omit the details.
\epf
\bigskip

Since  $I_b ,I_m $ and $g'(t), h'(t)$ are bounded in $(g(t),h(t))\times (0, T_0]$ by constants independent of $T_0$, then the local solution in $[0,T_0]$ to (\ref{a3}) can be extended for all $t\in(0,+\infty)$.

\begin{thm} Problem \eqref{a3} admits a global classical solution.
\end{thm}

Recalling that the system in (\ref{a3}) is quasimonotone nondecreasing, so the following comparison principle holds, see also Lemma 2.5 in \cite{ABL} or Lemma 3.5 in \cite{DL}.
\begin{lem} (The Comparison Principle)
  Suppose that $\overline g, \overline
h\in C^1([0, +\infty))$, $\overline I_b (x, t), \overline I_m (x,t)\in  C([\overline g(t), \overline h(t)]\times [0, +\infty))\cap
C^{2,1}((\overline g(t), \overline h(t))\times (0, +\infty))$, and
\begin{eqnarray*}
\left\{
\begin{array}{lll}
\frac{\partial \overline I_b }{\partial t}\geq D_1\frac{\partial^2 \overline I_b }{\partial x^2}-\gamma_{b}\overline I_b +\alpha_{b}\beta_{b}\frac{(N_b -\overline I_b)}{N_b} \overline I_m ,&\overline g(t)<x<\overline h(t),\, t>0, \\
\frac{\partial \overline I_m}{\partial t}\geq -d_{m}\overline I_m +\alpha_{m}\beta_{b} \frac{(A_m -\overline I_m )}{N_b} \overline I_b ,&\overline g(t)<x<\overline h(t),\, t>0,\\
\overline I_b (x,t)=\overline I_m (x, t)=0,\, & x=\overline g(t)\, \textrm{or}\, x= \overline h(t),\, t>0,\\
\overline g(0)\leq -h_0,\; \overline g'(t)\leq -\mu \frac{\partial \overline I_b }{\partial x}(\overline g(t), t), & t>0, \\
 \overline h(0)\geq h_0, \; \overline h'(t)\geq -\mu \frac{\partial \overline I_b }{\partial x}(\overline h(t), t), & t>0,\\
\overline I_b (x,0)\geq I_{b,0}(x),\ \overline I_m (x,0)\geq I_{m,0}(x),&-h_0\leq x\leq h_0.
\end{array} \right.
\end{eqnarray*}
Then the solution $(I_b , I_m ; g, h)$ of the free boundary problem $(\ref{a3})$ satisfies
$$h(t)\leq\overline h(t),\ g(t)\geq \overline g(t),\quad t\in [0, +\infty),$$
$$I_b (x, t)\leq \overline I_b (x, t),\ I_m (x, t)\leq \overline I_m (x, t),\quad x\in [g(t), h(t)],\, t\geq 0.$$\label{Com}
\end{lem}
\begin{rmk} The solution $(\overline I_b , \overline I_m ; \overline h, \overline g)$ in Lemma $\ref{Com}$ is usually called an upper solution
of \eqref{a3}. We can define a lower solution by
reversing all of the inequalities in the obvious places. Moreover, one
can easily prove an analogue of Lemma $\ref{Com}$ for lower solution.
\end{rmk}

 Next, we write $(I_b ^{\mu}, I_m ^{\mu}; g^{\mu}, h^{\mu})$ to examine the impact of $\mu$ on the solution, Lemma \ref{Com} leads directly to the following result.

\begin{cor} Let $I_{b,0}, I_{m,0}$ and other parameters and constants in $(\ref{a3})$ are fixed except $\mu$. If $\mu_1\leq \mu_2$, then $I_b ^{\mu_1}(x, t)\leq I_b ^{\mu_2}(x, t)$ and
$I_m ^{\mu_1}(x, t)\leq I_m ^{\mu_2}(x, t)$ over $\{(x,t):\, g^{\mu_1}(t)\leq x\leq h^{\mu_1}(t), t\geq 0\}$,
$g^{\mu_1}(t)\geq g^{\mu_2}(t)$ and $h^{\mu_1}(t)\leq h^{\mu_2}(t)$ in $[0, \infty)$.
\end{cor}

\section{The vanishing of WNv}
In this section, we concern about the conditions for vanishing of the virus.
 According to Lemma \ref{mono}, one can see that $x=h(t)$ is strictly increasing and $x=g(t)$ is strictly decreasing, therefore,
there exist $h_\infty, -g_\infty\in (0, +\infty]$ such that $\lim_{t\to +\infty} \ h(t)=h_\infty$
and $\lim_{t\to +\infty} \ g(t)=g_\infty$. We first present some properties of the free boundaries.
\begin{lem} Let $(I_b , I_m ; g, h)$ be a solution to $(\ref{a3})$
defined for $t\in[0, +\infty)$ and $x\in[g(t), h(t)]$. Then we have
$$-2h_0<g(t)+h(t)<2h_0 \mbox{ for } t\in[0, +\infty),$$
it means that the double moving fronts $x=g(t)$ and $x=h(t)$ are both finite or infinite simultaneously.
\end{lem}
\bpf It follows from continuity that $g(t)+h(t)>-2h_0$ holds for small $t>0$.
Let
$$T:=\sup\{s: g(t)+h(t)>-2h_0\ \mbox{ for all }\  t\in[0,s)\}.$$ As in \cite{ABL, DB}, we can assert
that $T= +\infty$. Otherwise, if $0<T< +\infty$ and
$$g(t)+h(t)>-2h_0 \mbox{ for } t\in[0,T),\ g(T)+h(T)=-2h_0.$$
We then have
\begin{eqnarray}
g'(T)+h'(T)\leq0. \label{Hq}
\end{eqnarray}

On the other hand, we define the functions
$$u(x, t):=I_b (x,t)-I_b (-x-2h_0, t),\ v(x, t):=I_m (x, t)-I_m (-x-2h_0, t)$$
over the region $$\Lambda:=\{(x, t):\ x\in[g(t), -h_0],\, t\in[0, T]\}.$$
It is easy to see that the pair $(u, v)$ is
well-defined for $(x, t)\in\Lambda$ since $ -h_0\leq -x-2h_0\leq -g(t)-2h_0\leq h(t) $, and the pair satisfies, for $g(t)<x<-h_0$, $0<t\leq T$,
$$u_t-d u_{xx}=-\gamma_{b} u+\alpha_{b}\beta_{b}\Big[\frac{(N_{b}-I_{b}(x,t))}{N_{b}} v-I_{m}(-x-2h_0,t)\frac{u}{N_b}\Big],$$
$$v_t=-d_{m} v+\alpha_{m}\beta_{b}\Big[\frac{(A_{m}-I_{m}(x,t))}{N_{b}} u-I_{b}(-x-2h_0,t)\frac{v}{N_b}\Big]$$
with
$$u(-h_0, t)=v(-h_0,t)=0,\, u(g(t), t)< 0,\, v(g(t), t)< 0 \mbox{ for } 0<t\leq T.$$ Moreover,
$$u(g(T), T)=I_b (g(T), T)-I_b (-g(T)-2h_0, T)=I_b (g(T), T)-I_b(h(T), T)=0.$$
Using the similar proof of Hopf boundary lemma, we get that
$$u(x,t)<0,\, v(x,t)<0  \mbox{ in } (g(t),-h_0)\times (0, T] \mbox{ and } u_x(g(T), T)<0.$$
Additionally, $$u_x(g(T), T)=\frac {\partial I_b }{\partial x}(g(T), T)+\frac {\partial I_b }{\partial x}(
h(T), T)=-[g'(T)+h'(T)]/\mu,$$
which implies that $$g'(T)+h'(T)>0,$$
we then leads a contradiction to (\ref{Hq}). So $T= +\infty$ and
$$g(t)+h(t)>-2h_0 \mbox{
for all } t>0.$$

Similarly, we can prove $g(t)+h(t)<2h_0$ for all $t>0$ by defining
$$u(x, t):=I_b (x,t)-I_b (2h_0-x, t),\ v(x, t):=I_m (x, t)-I_m (2h_0-x, t)$$ over the
region $\tilde \Lambda:=\{(x, t):\, x\in[h_0, h(t)],\, t\in[0, \tilde T]\}$ with
$\tilde T:=\sup\{s:\, g(t)+h(t)<2h_0 \mbox{ for all }  t\in[0,s)\}$. The proof is completed.  \epf
\bigskip

It follows from Lemma \ref{mono} that the infected habitat is expanding. Epidemically, if the infected habitat is limited and the infected cases disappear gradually,
we say the virus is vanishing and the epidemic is controlled. Mathematically, we have following definitions.
\begin{defi}
The virus is {\bf vanishing} if
$$h_\infty-g_\infty <\infty\ \textrm{ and}\
 \lim_{t\to +\infty} \ (||I_b (\cdot, t)||_{C([g(t),h(t)])}+||I_m (\cdot, t)||_{C([g(t), h(t)])})=0,$$
  and  {\bf spreading} if $$h_\infty-g_\infty =\infty\ \textrm{and}\
\limsup_{t\to +\infty}\ (||I_b (\cdot, t)||_{C([g(t),h(t)])}+||I_m (\cdot, t)||_{C([g(t),h(t)])})>0.$$
\end{defi}

 The following result shows that if $h_\infty-g_\infty<\infty$, then vanishing occurs.

 \begin{lem}\label{limit} If $h_\infty-g_\infty<\infty$, then there exists $\hat C$ independent of $t$ such that \label{vash}
 \begin{eqnarray}
\|I_b (\cdot, t)\|_{C^{1}([g(t), h(t)])}\leq \hat C,\ t\geq 1 ,\label{est-2}\\
||h'||_{C^{\alpha/2}([1, +\infty))},\, ||g'||_{C^{\alpha/2}([1, +\infty))}\leq \hat C.
\label{est-1}
\end{eqnarray}
Moreover,
 \begin{equation}
 \lim_{t\to
+\infty} \ (||I_b (\cdot, t)||_{C([g(t),h(t)])}+||I_m (\cdot, t)||_{C([g(t),h(t)])})=0.
\label{lim-1}
\end{equation}
\end{lem}
\bpf Similar to the proof of Theorem 2.1, we consider a transformation  $$y=\frac{2h_0x}{h(t)-g(t)}-\frac{h_0(h(t)+g(t))}{h(t)-g(t)},\
u(y,t)=I_b (x,t),$$ which straightens the free boundaries $x=h(t)$ and $x=g(t)$ to fixed lines $y=h_0$ and $y=-h_0$ respectively.
Hence the free boundary problem (\ref{a3}) becomes the fixed boundary problem (\ref{Hb}).
Since $-g(t)$ and $h(t)$ are increasing and bounded, it follows from the standard $L^p$ theory and the Sobolev imbedding
theorem (\cite{LSU, Lie}) that for $0<\alpha <1$,
there exists a constant $\tilde C$
depending on $\alpha, h_0, \|I_{b,0}\|_{C^{2}[-h_0, h_0]}$, $\|I_{m,0}\|_{C^{2}[-h_0, h_0]}$, and $g_\infty, h_\infty$ such that
\begin{eqnarray}\|u\|_{C^{1+\alpha,
(1+\alpha)/2}([-h_0, h_0]\times [\tau, \tau+1])}\leq \tilde C\label{Bg1}
\end{eqnarray}
for any $\tau\geq 1$. Recalling that $\tilde C$ is independent of $\tau$ and $g'(t), h'(t)$ are bounded by $C_1$ from Lemma 2.2, we then arrive at (\ref{est-2}) and (\ref{est-1}).
Using (\ref{est-1}) and the assumption that $h_\infty-g_\infty<\infty$ gives
$$h'(t)\to 0\ \textrm{and}\ g'(t)\to 0\ \textrm{as}\, t\to +\infty.$$

Next we are going to derive (\ref{lim-1}). Suppose that
$$\limsup_{t\to +\infty} \ ||I_b (\cdot, t)||_{C([g(t), h(t)])}=\delta>0$$
 by contradiction. Then there exists a sequence $\{(x_k, t_k )\}$
in $(g(t), h(t))\times (0, \infty)$
such that $I_b (x_k,t_k)\geq \delta /2$ for all $k \in \mathbb{N}$, and $t_k\to \infty$ as $k\to \infty$.
Since that  $-\infty<g_\infty<g(t)<x_k<h(t)<h_\infty<\infty$, we can extract a subsequence of $\{x_k\}$ (still denoted by it),
such that $x_k\to x_0\in [g_\infty, h_\infty]$ as $k\to \infty$.

Due to the uniform boundedness in (\ref{est-2}), we assert that $x_0\in (g_\infty, h_\infty)$. In fact, if $x_0=h_\infty$, then
$(x_k-h(t_k))\to (x_0-h_\infty)=0$ as $k\to \infty$. On the other hand,
$$ \delta /2\leq I_b (x_k, t_k)=I_b (x_k, t_k)-I_b (h(t_k),t_k)=\frac{\partial I_b}{\partial x}(\xi_k, t_k)(x_k-h(t_k))\leq -\hat C (x_k-h(t_k)),$$
where $\xi_k\in (x_k, h(t_k))$. So, $(h(t_k)-x_k)\geq \frac {\delta}{2C_1}$ for $k\in \mathbb{N}$, which leads to a contradiction and then $x_0\neq h_\infty$.
Similarly, we have $x_0\neq g_\infty$.

Let $u_k(x,t)=I_b (x,t_k+t)$ and $v_k(x,t)=I_m (x,t_k+t)$  for
$x\in (g(t_k+t), h(t_k+t)), t\in (-t_k, \infty)$.
By the parabolic regularity, we deduce that  $\{(u_k, v_k)\}$ has a subsequence $\{(u_{k_i}, v_{k_i})\}$ such that
$(u_{k_i}, v_{k_i})\to (\tilde u, \tilde v)$ as $i\to \infty$ and $(\tilde u, \tilde v)$ satisfies
\begin{eqnarray*} \left\{
\begin{array}{lll}
\tilde u_t-D_1 \tilde u_{xx}=-\gamma_{b}\tilde u+\alpha_{b}\beta_{b}\frac{(N_b - \tilde u)}{N_b} \tilde v,\; & g_\infty<x<h_\infty,\ t\in (-\infty, \infty),  \\
\tilde v_t=-d_{m}\tilde v+\alpha_{m}\beta_{b} \frac{(A_{m} -\tilde v)}{N_b} \tilde u,\; &\ g_\infty<x<h_\infty, \ t\in (-\infty, \infty).
\end{array} \right.
\end{eqnarray*}
Recalling that $u_{k}(x_k,0)=I_b(x_k, t_k)\geq \delta/2$ for any $k$, we then have $\tilde u(x_0, 0)\geq \delta/2$.
Furthermore, $\tilde u>0$ in $ (g_\infty, h_\infty)\times(-\infty, \infty)$.
 Applying the Hopf lemma at the point $(h_\infty, 0)$ yields
 $\tilde u_x(h_\infty, 0 )\leq -\sigma_0$ for some $\sigma_0>0$.
On the other hand, $h'(t)\to 0$ as $t\to \infty$, that is,
$\frac {\partial I_b }{\partial x}(h(t_k),t_k)\to 0$ as $t_k\to \infty$ by the free boundary condition. Using (\ref{Bg1}), which suggests that $I_b $ has a uniform
$C^{1+\alpha, (1+\alpha)/2}$ bound over $\{(x,t): g(t)\leq x\leq h(t),\, t\geq 1\}$, we then derive
$\frac {\partial I_b }{\partial x}(h(t_k),t_k+0)=(u_k)_x(h(t_k),0)\to \tilde u_x(h_\infty,0)$ as $k\to \infty$, and therefore $\tilde u_x(h_\infty,0)=0$,
which contradicts the fact that  $\tilde u_x(h_\infty,0)\leq -\sigma_0<0$.
 Thus $\lim_{t\to +\infty} \ ||I_b (\cdot,t)||_{C([g(t),h(t)])}=0$.

Noting that $I_m (x,t)$ satisfies
 $$\frac{\partial I_m(x,t)}{\partial t}= -d_{m}I_m (x,t)+\alpha_{m}\beta_{b} \frac{(A_m -I_m (x,t))}{N_b} I_b (x,t),\ g(t)<x<h(t),\, t>0,$$
 and $\alpha_{m}\beta_{b} \frac{(A_m -I_m (x,t))}{N_b} I_b (x,t)\to 0$ uniformly for $x\in [g(t), h(t)]$ as $t\to \infty$, we immediately have that $\lim_{t\to
+\infty} \ ||I_m (\cdot, t)||_{C([g(t),h(t)])}=0$.
\epf

\bigskip
As noted in the introduction section, when we consider the spreading or vanishing of the virus, a threshold parameter $R_0$ is usually
defined for differential systems describing epidemic models. $R_0$ is called the basic reproduction number. But in our model \eqref{a3}, the infected interval is changing with the time $t$, therefore,
the basic reproduction number is not a constant and
should be a function of $t$. So we here call it the risk index, which is expressed by
\begin{equation}\label{free}
R_0^F(t):=R_0^D((g(t), h(t)))=\frac {{\alpha_{m}\alpha_{b} \beta_{b}^2  A_m}}{{N_{b}}{d_{m}}({\gamma_{b}+D_1(\frac \pi {h(t)-g(t)})^2})},
\end{equation}
where $R_0 ^D ( \Omega)$ is a threshold parameter for the corresponding problem (\ref{Aa1}) in $\Omega$
with null Dirichlet boundary condition on $\partial\Omega$, and it depends on the principal eigenvalue of the corresponding problem.
With the above definition, we have the following properties of $R_0^F (t)$.
\begin{lem} The following statements are valid:\label{proper}
\item[$(i)$] ${\textrm sgn}(1-R_0^F(t))$ = ${\textrm sgn} \lambda_0 $, where $\lambda_0$ is the principal eigenvalue of the problem
\begin{eqnarray}
\left\{
\begin{array}{lll}
-D_1 \psi_{xx}=-\gamma_{b} \psi+\frac {\alpha_{m}\alpha_{b}\beta_{b}^2 A_m}{N_{b} d_m}\psi+\lambda_0 \psi,\; &
x\in (g(t), h(t)),  \\
\psi(x)=0, &x=g(t)\ \textrm{or}\ x=h(t);
\end{array} \right.
\label{B1f}
\end{eqnarray}

\item[$(ii)$] $R_0^F(t)$ is strictly monotone increasing function of $t$, that is if $t_1<t_2$, then $R_0^F(t_1)<R_0^F(t_2)$;

\item[$(iii)$] if $h(t)-g(t)\to \infty$ as $t\to \infty$, then $R_0^F(t)\to R_0$ as $t\to \infty$.
\end{lem}
\bpf $(ii)$ and $(iii)$ can be obtained directly from the expression (\ref{free}). As to $(i)$,
direct calculation shows that
$$\lambda_0=\gamma_{b}+D_1(\frac \pi {h(t)-g(t)})^2-\frac {\alpha_{m}\alpha_{b}\beta_{b}^2 A_m}{N_{b} d_m}$$\\
$$=\big[\gamma_{b}+D_1(\frac \pi {h(t)-g(t)})^2\big](1-R_0^F(t)),$$
which implies that $(i)$ holds.
\epf
\bigskip

In the following, we will explore some effective ways to control the virus. Mathematically, we discuss sufficient conditions so that the virus is vanishing.
\begin{thm} If $R_0\leq 1$, then $h_\infty-g_\infty<\infty$ and
$\lim_{t\to +\infty} \ (||I_b (\cdot, t)||_{C([g(t),h(t)])}+||I_m (\cdot, t)||_{C([g(t),h(t)])})=0$.
\label{vanish}
\end{thm}
\bpf
We first prove that $h_\infty -g_\infty<+\infty$. In fact, direct computations yield
\begin{eqnarray*}& &\frac{\textrm{d}}{\textrm{d} t}\int_{g(t)}^{h(t)}\Big[I_b (x, t)+\frac {\alpha_{b}\beta_{b}}{d_m}  I_m (x, t)\Big]\textrm{d}x\\
&=&\int_{g(t)}^{h(t)}\Big[\frac{\partial I_b}{\partial t}+\frac {\alpha_{b}\beta_{b}}{d_m}\frac{\partial I_m}{\partial t}\Big](x, t)\textrm{d}x+h'(t)\Big[I_b +\frac {\alpha_{b}\beta_{b}}{d_m} I_m \Big](h(t), t)\\
& &\quad
-g'(t)\Big[I_b+\frac {\alpha_{b}\beta_{b}}{d_m} I_m\Big](g(t), t)\\[1mm]
&\leq&\int_{g(t)}^{h(t)}D_1 \frac{\partial^2 I_b}{\partial x^2}\textrm{d}x+\int_{g(t)}^{h(t)}\Big[-\gamma_{b} I_b (x,t)+\frac {\alpha_{m}\alpha_{b} \beta_{b}^2 A_m}{N_{b} d_m} I_b (x, t)\Big]\textrm{d}x\\[1mm]
&\leq&-\frac{D_1}{\mu}(h'(t)-g'(t))+\int_{g(t)}^{h(t)}\Big[-\gamma_{b}I_b (x,t)+\frac {\alpha_{m} \alpha_{b} \beta_{b}^2 A_m}{N_{b} d_m} I_b (x, t)\Big] \textrm{d}x.
\end{eqnarray*}
Integrating from $0$ to $t\,(>0)$ gives
\begin{eqnarray}
& &\int_{g(t)}^{h(t)}\Big[I_b+\frac {\alpha_{b} \beta_{b}}{d_{m}}I_m\Big](x, t)\textrm{d}x \nonumber\\
&\leq& \int ^{h(0)}_{g(0)}\Big[I_b+\frac {\alpha_{m} \beta_{b}}{d_{m}}I_m\Big](x, 0)\textrm{d}x+\frac {D_1}{\mu}(h(0)-g(0))-\frac {D_1}{\mu}(h(t)-g(t))\nonumber \\
& &+\int_{0}^t\int_{g(s)}^{h(s)}\Big[-\gamma_{b} I_b (x,t)+\frac {\alpha_{m} \alpha_{b} \beta_{b}^2 A_{m}}{N_{b} d_{m}} I_b (x, t)\Big]dxds,
\quad t\geq 0.\label{k1}
\end{eqnarray}
It follows from $R_0\leq 1$ that $-\gamma_{b} I_b (x,t)+\frac {\alpha_{m}\alpha_{b} \beta_{b}^2 A_m}{N_b d_m}I_b (x, t)\leq 0$ for $x\in [g(t),h(t)]$ and $t\geq 0$, we then have
$$\frac {D_1}{\mu}(h(t)-g(t)) \leq \int ^{h(0)}_{g(0)}\big[I_b +\frac {\alpha_{b} \beta_{b}} {d_m}  I_m\big](x, 0)\textrm{d}x
+\frac {D_1}{\mu}(h(0)-g(0))$$
for $t\geq 0$, which implies that $h_\infty-g_\infty<\infty$. Furthermore, the vanishing of the virus
 follows easily from Lemma 3.2.
\epf

\begin{thm} \label{small}If $R_0^F(0)<1$ and $\mu$ is sufficiently small. Then $h_\infty-g_\infty<\infty$ and
$\lim_{t\to +\infty} \ (||I_{b} (\cdot, t)||_{C([g(t),h(t)])}+||I_{m} (\cdot, t)||_{C([g(t),h(t)])})=0$.
\end{thm}
\bpf We are going to  construct a suitable upper solution to problem \eqref{a3}.
Since $R_0^F(0)<1$, it follows from Lemma \ref{proper} that there exist $\lambda_0>0$ and $0<\psi(x)\leq 1$ in $(-h_0, h_0)$ such that
\begin{eqnarray}
\left\{
\begin{array}{lll}
-D_1 \psi_{xx}= -\gamma_{b}\psi+\frac {\alpha_{m} \alpha_{b} \beta_{b}^2 A_m}{N_{b} d_m}\psi+\lambda_0 \psi,\; &
-h_0<x<h_0,  \\
\psi(x)=0, &x=\pm h_0.
\end{array} \right.
\label{B1f1}
\end{eqnarray}
Accordingly, there exists a small $\delta_1 >0$ such that, for $\delta\leq \delta_1$,
$$-\delta +(\frac 1{(1+\delta)^2}-1)\frac {\alpha_{m} \alpha_{b} \beta_{b}^2 A_m}{N_{b} d_m}+[\frac 1{(1+\delta)^2}-\frac {1}{ 2 }]\lambda_0\geq 0.$$
Similarly as in \cite{ABL}, we set
$$\sigma (t)=h_0(1+\delta-\frac \delta 2 e^{-\delta t}), \  t\geq 0,$$
and
$$\overline I_b (x, t)=M e^{-\delta t}\psi(xh_0/\sigma (t)), \ -\sigma(t)\leq
x\leq \sigma(t),\ t\geq 0.$$
$$\overline I_m (x,t)=(\frac {\alpha_{m} \beta_{b} A_m} {N_b d_m}+\frac{\lambda_0}{2\alpha_b\beta_b})\overline I_b (x,t), \ -\sigma(t)\leq
x\leq \sigma(t),\ t\geq 0.$$
Direct calculations give
\begin{eqnarray*}
& &\dfrac{\partial \overline I_b}{\partial t}-D_1 \dfrac{\partial^2\overline I_b}{\partial x^2}+
\gamma_{b}\overline I_b -\alpha_{b} \beta_{b}\frac{(N_b-\overline I_b)}{N_b} \overline I_m \\
&\geq &\dfrac{\partial \overline I_b}{\partial t}-D_1 \dfrac{\partial^2\overline I_b}{\partial x^2}+
\gamma_{b}\overline I_b -\alpha_{b} \beta_{b} \overline I_m \\
& =& -\delta \overline I_b-M e^{-\delta t}\psi'\frac{xh_0\sigma'(t)}{\sigma^2(t)}+\Big(\frac{h_0}{\sigma(t)}\Big)^2
\Big[-\gamma_{b}+\frac {\alpha_{m} \alpha_{b} \beta_{b}^2 A_m}{N_{b}d_m}+\lambda_0 \Big]\overline I_b \\
& & +\Big[\gamma_{b}-\frac {\alpha_{m} \alpha_{b} \beta_{b}^2 A_m}{N_{b} d_m}-\frac{\lambda_0}{4 }\Big]\overline I_b \\
& \geq& \overline I_b \Big\{-\delta +(\frac 1{(1+\delta)^2}-1)\frac {\alpha_{m} \alpha_{b} \beta_{b}^2 A_m}{N_{b}d_m}+\big[\frac 1{(1+\delta)^2}-\frac{1}{2}\big]\lambda_0\Big\}\geq 0,
\end{eqnarray*}
\begin{eqnarray*}
& &\dfrac{\partial \overline I_m}{\partial t}+d_{m}\overline I_m-\alpha_{m} \beta_{b} \frac{(A_{m}-\overline I_m)}{N_b}\overline I_b \\
& \geq&\dfrac{\partial \overline I_m}{\partial t}+d_{m}\overline I_m-\frac{\alpha_{m} \beta_{b} A_{m}}{N_{b}} \overline I_b \\
& =& -\delta \overline I_m-M e^{-\delta t}\psi'\frac{xh_0\sigma'(t)}{\sigma^2(t)}\Big(\frac {\alpha_{m} \beta_{b} A_m}{N_{b} d_m}+\frac{\lambda_0}{2\alpha_{b} \beta_{b}}\Big)\\
& &+d_{m}\Big(\frac {\alpha_{m} \beta_{b} A_{m}}{N_b d_m}+\frac{\lambda_0}{2\alpha_{b} \beta_{b}}\Big) \overline I_b -\frac{\alpha_{m} \beta_{b} A_m}{N_b}\overline I_b \\
& \geq& (d_{m}-\delta )\Big(\frac {\alpha_{m} \beta_{b} A_m}{N_{b}d_m}+\frac{\lambda_0}{2\alpha_{b} \beta_{b}}\Big)\overline I_b -\frac{\alpha_{m} \beta_{b} A_m}{N_b}\overline I_b \\
& =&\overline I_b\Big\{{d_m \frac{\lambda_0}{2\alpha_{b} \beta_{b}}-\delta\big[\frac{\alpha_m \beta_b A_m}{N_b d_m}+\frac{\lambda_0}{2\alpha_{b} \beta_{b}}\big]}\Big\}
\end{eqnarray*}
for all $-\sigma (t)<x<\sigma (t)$ and $t>0$. So, taking
$$\delta=\min\{\delta_1,\,  \frac{d_m\lambda_0}{2\alpha_{b} \beta_{b}}[\frac{\alpha_m \beta_b A_m}{N_b d_m}+\frac{\lambda_0}{2\alpha_{b} \beta_{b}}]^{-1}\},$$
we then have
\begin{eqnarray*}
\left\{
\begin{array}{lll}
\dfrac{\partial \overline I_b}{\partial t}\geq D_1 \dfrac{\partial^2 \overline I_b}{\partial x^2}
-\gamma_{b}\overline I_b+\alpha_{b} \beta_{b}\frac{(N_b -\overline I_b)}{N_b} \overline I_m,\; & -\sigma(t)<x<\sigma(t),\, t>0, \\
\dfrac{\partial \overline I_m}{\partial t}\geq -d_{m}\overline I_m +\alpha_{m} \beta_{b} \frac{(A_m -\overline I_m)}{N_b} \overline I_b,\; &  -\sigma(t)<x<\sigma(t), \, t>0,\\
\overline I_b (x,t)=\overline I_m (x, t)=0,&x=\pm \sigma(t),\, \, t>0.
\end{array} \right.
\end{eqnarray*}
Now, we can choose $M$ big enough so that $I_{b,0}(x)\leq M \psi(\frac {h_0}{1+\delta/2})\leq \overline I_b (x, 0)=M \psi(\frac {x}{1+\delta/2})$
and $I_{m,0}(x)\leq M \psi(\frac {h_0}{1+\delta/2})(\frac {\alpha_{m} \beta_{b} A_m}{N_{b} d_m}+\frac{\lambda_0}{2\alpha_{b} \beta_{b}})\leq \overline I_m (x, 0)$
 for $x\in [-h_0, h_0]$.

Additionally, since that
$$\sigma'(t)=h_0 \frac {\delta^2} 2 e^{-\delta t},\quad -\frac{\partial\overline I_b}{\partial x}(
\sigma (t),t)=-M \frac {h_0}{\sigma (t)}\psi'(h_0)e^{-\delta t},$$
$$-\frac{\partial\overline I_b}{\partial x}(
-\sigma (t),t)=-M \frac {h_0}{\sigma (t)}\psi'(-h_0)e^{-\delta t},\quad \psi'(-h_0)=-\psi'(h_0),$$
we can choose $\mu=-\frac {\delta^2h_0(1+\delta)} {2M \psi'(h_0)}$ such that
$$ -\sigma'(t)\leq -\mu \frac{\partial \overline I_b}{\partial x}(-\sigma(t), t),\quad
\sigma'(t)\geq -\mu \frac{\partial \overline I_b}{\partial x}(\sigma(t), t)$$
for $t>0$. Hence, Using Lemma 2.4 concludes that $g(t)\geq -\sigma(t)$ and $h(t)\leq\sigma(t)$ for $t>0$. It
follows that $h_\infty-g_\infty\leq \lim_{t\to\infty}
2\sigma(t)=2h_0(1+\delta)<\infty$, and then $\lim_{t\to +\infty} \ (||I_b (\cdot, t)||_{C([g(t),h(t)])}+||I_m (\cdot, t)||_{C([g(t),h(t)])})=0$ by Lemma \ref{vash}.
 \epf

 \bigskip

 Using the similar upper solution, we can also prove that vanishing happens for small initial data.
\begin{thm} If $R_0^F(0)<1$ and the initial functions $I_{b,0}(x)$ and $I_{m,0}(x)$ are sufficiently small. Then $h_\infty-g_\infty<\infty$ and
$\lim_{t\to +\infty} \ (||I_{b} (\cdot, t)||_{C([g(t),h(t)])}+||I_{m} (\cdot, t)||_{C([g(t),h(t)])})=0$.
\end{thm}

\section{The spreading of WNv}

In this section, our aim is to look for some factors which lead to the spreading of the virus.
\begin{thm} If $R_0^F(t_0)\geq 1$ for $t_0\geq 0$, then
$h_\infty=-g_\infty=\infty$ and
$$\liminf_{t\to
+\infty} \ ||I_b (\cdot, t)||_{C([g(t), h(t)])}>0,$$
 that is, spreading must occur.
\label{spread}
\end{thm}
\bpf It suffices to prove it in the case $R_0 ^F (t_0)>1$. Since if $R_0^F(t_0)=1$, for any given $t_1>t_0$, we then have $g(t_1)<g(t_0)$ and $h(t_1)>h(t_0)$, which yields $R_0^F(t_1)>
R_0^F(t_0)=1$ from the monotonicity in Lemma \ref{proper}.
Hence replacing the time $t_0$ by $t_1$, we can obtain $h_\infty-g_\infty=+\infty$ as the following.

In this case $R_0 ^F (t_0)>1$, the following eigenvalue problem
\begin{eqnarray}
\left\{
\begin{array}{lll}
- D_1\psi_{xx}= -\gamma_{b}\psi+\frac {\alpha_{m} \alpha_{b} \beta_{b}^2 A_m}{N_{b} d_m}\psi+\lambda_0 \psi,\; &
g(t_0)<x<h(t_0),  \\
\psi(x)=0, &x=g(t_0)\ \textrm{or}\ x=h(t_0)
\end{array} \right.
\label{B2f}
\end{eqnarray}
 admits a positive solution $\psi(x)$ with $||\psi||_{L^\infty}=1$, and the principal eigenvalue $\lambda_0<0$ by Lemma \ref{proper}.

Next, we are going to construct a suitable lower solution to
\eqref{a3}, and define
$$\underline I_{b}(x,t)=\delta \psi(x),\quad \underline I_m=\Big(\frac {\alpha_{m} \beta_{b} A_m}{N_{b} d_m}+\frac{\lambda_0}{2 \alpha_{b}\beta_{b}}\Big)\delta \psi(x)$$
for $g(t_0)\leq x\leq h(t_0)$, $t\geq t_0$, where $\delta $ is chosen later.

 It follows from the direct calculations that
\begin{eqnarray*}
& & \dfrac{\partial \underline I_b}{\partial t}-D_1 \dfrac{\partial^2\underline I_b}{\partial x^2}+\gamma_{b}
\underline I_b-\alpha_{b} \beta_{b}\frac{(N_b-\underline I_b)}{N_b} \underline I_m \\
&=&\delta\psi(x) \left\{\frac{1}{2} \lambda_0 +\delta\psi\frac{\alpha_{b}\beta_{b}}{N_b} \big[\frac{\alpha_{m} \beta_{b}A_m}{N_b d_m}+\frac{\lambda_0}{2\alpha_b \beta_b}\big] \right\},\\
& &\dfrac{\partial \underline I_m}{\partial t}+d_{m}\underline I_m-\alpha_{m} \beta_{b} \frac{(A_{m}- \underline I_m)}{N_b} \underline I_b\\
&=& \delta\psi(x)\left\{\frac{d_m}{2\alpha_b \beta_b} \lambda_0 +\delta\psi\frac{\alpha_{m}\beta_{b}}{N_b} \big[\frac{\alpha_{m} \beta_{b}A_m}{N_b d_m}+\frac{\lambda_0}{2\alpha_b \beta_b}\big] \right\}
\end{eqnarray*}
for all $g(t_0)<x<h(t_0)$ and $t>t_0$. Noting that $\lambda_0<0$ and $0\leq \psi(x)\leq 1$, we can chose $\delta $ sufficiently small such that
\begin{eqnarray*}
\left\{
\begin{array}{lll}
\dfrac{\partial \underline I_b}{\partial t}\leq D_1 \dfrac{\partial^2 \underline I_b}{\partial x^2}
-\gamma_{b}\underline I_b +\alpha_{b} \beta_{b}\frac{(N_b -\underline I_b)}{N_b}\underline I_m,\; & g(t_0)<x<h(t_0),\, t>t_0, \\
\dfrac{\partial \underline I_m}{\partial t}\leq -d_{m}\underline I_m +\alpha_{m} \beta_{b} \frac{(A_m -\underline I_m)}{N_b} \underline I_b ,\; &  g(t_0)<x<h(t_0), \, t>t_0,\\
\underline I_b (x,t)=\underline I_m (x, t)=0,&x=g(t_0)\ \textrm{or}\, x=h(t_0)\, \, t>t_0,\\
0=(g(t_0))'\geq -\mu \frac{\partial \underline I_b}{\partial x}(g(t_0), t), & t>t_0, \\
0=(h(t_0))'\leq -\mu \frac{\partial \underline I_b}{\partial x}(h(t_0), t), & t>t_0,\\
\underline {I_b}(x, t_0)\leq I_{b,0}(x),\ \underline{I_m}(x,t_0)\leq I_{m,0}(x),\; &g(t_0)\leq x\leq h(t_0).
\end{array} \right.
\end{eqnarray*}
Thus, using Remark 2.1 gives that $I_b (x,t)\geq\underline I_b (x,t)$ and  $I_m (x,t)\geq\underline I_m (x,t)$
in $[g(t_0), h(t_0)]\times [t_0, \infty)$. It follows that $\liminf_{t\to
+\infty} \ ||I_b (\cdot, t)||_{C([g(t), h(t)])}\geq \delta \psi(0)>0$ and then $h_\infty-g_\infty=+\infty$ by Lemma \ref{vash}.
 \epf

\begin{rmk} Theorem $\ref{vanish}$ shows that vanishing always happens for $R_0\leq 1$.
If $R_0>1$, $R_0^F(t_0)\geq 1$ is equivalent to that $(h(t_0)-g(t_0))\geq \pi \sqrt{\frac{D_1}{\gamma_{b}(R_0 -1)}}$.
 Theorem $\ref{spread}$ reveals a critical spreading length, which may be called a
``spreading barrier'', $l^*=\pi\sqrt{\frac{D_1}{\gamma_{b}(R_0 -1)}}$,
 such that the virus will spreads to all the new
population if its spreading length can break through this barrier $l^*$ in some finite time, or
the spreading never breaks through this barrier and the virus vanishes in the long run.
\end{rmk}

Recalling Theorem 3.5, we know that a small expanding rate $\mu$ is benefit for the vanishing of the virus. We wonder what will happen for the
 virus if $\mu$ becomes large. For this purpose, we first consider the following initial boundary value problem
\begin{eqnarray}
\left\{
\begin{array}{lll}
u_{t}-D_1 u_{xx}=f(x,t)u,\; &g(t)<x<h(t),\; t>0,   \\
u(x,t)=0,&x=g(t)\, \textrm{or}\, x=h(t),\, t\geq0,\\
g'(t)=-\mu u_{x}(g(t),t),\;  & t>0, \\
h'(t)=-\mu u_{x}(h(t),t),\;  & t>0, \\
u(x,0)=u_{0}(x), & -h_0\leq x\leq h_0,
\end{array} \right.
\label{a991}
\end{eqnarray}
where $f(x,t)$ is a continuous function, $u_0\in C^2 [-h_0,h_0]$, $u_0(\pm h_0)=0$ and $u_0(x)>0, x\in (-h_0,h_0)$.
\begin{lem} \label{prop2} Assume that there exists a constant $M_{1}$ such that $f(x,t)\geq -M_{1} $ for $-\infty<x<\infty$, $t>0$.
 Then for any given constant $H>0$, there exists $\mu_H>0$, such that when $\mu>\mu_H$, the corresponding unique solution $(u^\mu (x,t); g^\mu (t)$, $h^\mu (t))$ of problem $($\ref{a991}$)$ satisfies
\begin{equation}\label{133}
\limsup_{t\to +\infty} \ g^\mu (t)<-H \quad and \quad
\liminf_{t\to +\infty} \ h^\mu (t)>H.
\end{equation}
\end{lem}
\bpf
   We start with the following initial-boundary value problem
 \begin{eqnarray}
\left\{
\begin{array}{lll}
v_{t}-D_1 v_{xx}=-M_1 v,\; &p(t)<x<q(t),\; t>0,   \\
v(x,t)=0,&x=p(t)\, \textrm{or}\, x=q(t),\, t\geq0,\\
p'(t)=-\mu v_{x}(p(t),t), p(0)=-h_0<0, \;  & t>0, \\
q'(t)=-\mu v_{x}(q(t),t), q(0)=h_0>0, \;  & t>0, \\
v(x,0)=u_{0}(x), & -h_0\leq x\leq h_0,
\end{array} \right.
\label{a992}
\end{eqnarray}
it admits a unique global solution $(v^\mu ; g^\mu , h^\mu)$ and $(p^\mu)'(t)<0$, $(q^\mu)'(t)<0$ for $t>0$.
It follows from Corollary 2.5 and the comparison principle that
\begin{equation}\label{134}
\begin{split}
&u^\mu(x,t)\geq v^\mu (x,t)  \quad \textrm{for}\       p^\mu (t)\leq x\leq q^\mu (t),\, t>0 \\
&g^\mu(t)\leq p^\mu (t),\quad h^\mu(t)\geq q^\mu (t)\quad \textrm{for}\   t>0.
\end{split}
\end{equation}

Now we are going to prove that for all large $\mu$,
 \begin{equation}\label{135}
p^\mu (2)\leq-H  \quad and\;  \quad q^\mu (2)\geq H.
\end{equation}
Choosing smooth functions $\underline p (t)$ and $\underline q (t)$ with
$$\underline p (0)=-\frac{h_0}{2},\quad \underline p (2)=-H$$
and
$$\underline q (0)=\frac{h_0}{2},\quad \underline q (2)=H,$$
we then consider the following problem
 \begin{eqnarray}
\left\{
\begin{array}{lll}
\underline v_{t} -D_1 v_{xx}=-M_1\underline v ,\; &\underline p (t) <x<\underline q(t) ,\; t>0,   \\
\underline v (x,t)=0,&x=\underline p (t)\, \textrm{or}\, x=\underline q (t),\, t\geq0,\\
\underline v (x,0)=\underline v_{0} (x), & -\frac{h_0}{2}\leq x\leq \frac{h_0}{2},
\end{array} \right.
\label{a993}
\end{eqnarray}
where the smooth value $\underline v_{0} (x)$ satisfies
\begin{eqnarray}
\left\{
\begin{array}{lll}
0<\underline v_{0} (x)<u_{0}(x), \; -\frac{h_0}{2}\leq x\leq \frac{h_0}{2},\\
\underline v_{0} (-\frac{h}{2})=v_{0}(\frac{h}{2})=0,\ \underline v'_{0} (-\frac{h}{2})>0,\ v'_{0}(\frac{h}{2})<0.
\end{array} \right.
\label{a994}
\end{eqnarray}
Hence, the standard theory for parabolic equations ensures that problem (\ref{a993}) has a unique solution $(\underline v;\underline p,\underline q)$ with ${\underline v}_x (\underline p (t),t)>0$ and ${\underline v}_x (\underline q (t),t)<0$ for $t\in[0,2]$ by using Hopf boundary lemma.

According to our choice of $\underline v _0 (x),\underline p (t)$ and $\underline q (t),$ there exists a constant $\mu_H,$ such that for all $\mu>\mu_H$,
\begin{equation}\label{136}
\underline p' (t)\geq -\mu v_{x}(\underline p (t),t) \;  and  \quad \; \underline q' (t)\leq -\mu v_{x}(\underline q (t),t), \; 0\leq t\leq 2.
\end{equation}
It is easy to see that,\\
$$\underline p (0)=-\frac{h}{2}>-h_0 =p^\mu (0), \; \underline q (0)=\frac{h}{2}<h_0 =q^\mu (0).$$\\
Using (\ref{a992}),(\ref{a993}),(\ref{a994}), (\ref{136}) and the comparison principle gives
$$v ^\mu(x,t)\geq \underline v (x,t), \quad p ^\mu(t)\leq \underline p (t) \quad and \quad q^\mu(t)\geq \underline q (t),$$ for $\underline p (t)\leq x\leq \underline q (t), 0\leq t\leq 2,$
which means that (\ref{135}) holds. Thanks to (\ref{134}) and (\ref{135}), we obtain
$$\limsup_{t\to +\infty} \ g (t)\leq \lim_{t\to +\infty} \ p^\mu (t)\leq p^\mu (2)\leq-H,$$
$$\liminf_{t\to +\infty} \ q (t)\geq \lim_{t\to +\infty} \ q^\mu (t)\geq q^\mu (2)\geq H.$$
\epf

\bigskip

 The following result
shows that spreading happens for large expanding capicity.
 \begin{thm} Assume that $R_0^F(0)<1<R_0$. Then $h_\infty-g_\infty=\infty$ and spreading happens for large $\mu$.
\end{thm}
\bpf
Recalling that $R^D_0((-L, L))\to R_0>1$ as $L\to +\infty$, then there exists $H>0$ such that
$R^D_0((-H, H))>1$.
For given $H$, since that
\begin{equation}\label{137}
\frac{\partial I_{b}}{\partial t}-D_1\frac{\partial^2I_b}{\partial x^2}\geq -\gamma_b I_{b}
\end{equation}
from the first equation in $(\ref{a3})$, using Lemma $4.2$ yields that there exists $\mu_H>0$ such that for any $\mu>\mu_H$,
\begin{equation}\label{138}
\limsup_{t\to +\infty} \ g (t)<-H  \ \textrm{and} \ \liminf_{t\to +\infty} \ h(t)>H,
\end{equation}
which together with the monotonicity of $g(t)$ and $h(t)$ gives that there exists $T_0>0$ such that $g(T_0)<-H $ and $h(T_0)>H$, therefore, we have
$$R_0 ^F (T_0)=R_0 ^D ((g(T_0),h(T_0)))>R_0 ^D ((-H,H))>1.$$
Thus, for large $\mu$, we can apply Theorem \ref{spread} to conclude that $h_\infty-g_\infty=\infty$ and the spreading happens.
\epf
\bigskip

Considering $\mu$ as a varying parameter, we have the next theorem.
\begin{thm} (Sharp threshold) For any fixed $h_0$, $I_{b,0}$ and $I_{m,0}$ satisfying $(\ref{Ae1})$, there exists $\mu^*\in [0, \infty]$
 such that spreading occurs when $\mu> \mu^*$, and vanishing occurs when $0<\mu\leq \mu^*$.
\end{thm}
\bpf
Theorem \ref{spread} shows that spreading always happens if $R_0^F(0)\geq 1$. Thus, in this
case we have $\mu^*=0$. Theorem \ref{vanish} shows that vanishing happens if $R_0\leq 1$, so in this case $\mu^*=\infty$.

For the remaining case $R_0^F(0)<1<R_0$, define
$$\sum=\{\mu>0: h_\infty^{\mu}-g_\infty^{\mu} \leq l^* \}\ \textrm{and}\ \mu^*:=\sup \sum,$$
where $l^*=\pi\sqrt{\frac{D_1}{\gamma_{b}(R_0 -1)}}$ defined in Remark 4.1.
Thanks to the monotonicity of $h(t)$ and $g(t)$ with respect to $\mu$ (Corollary 2.5), we see from Theorem 3.5 that the set $\sum$ is not empty and $\mu^*>0$,
it also follows from Theorem 4.3 that $\mu^*<\infty$. Therefore, the virus spreads if $\mu>\mu^*$ and vanishes if $0<\mu<\mu^*.$

We now claim that the vanishing happens for $\mu=\mu^*$.
Otherwise $h^{\mu^*}_\infty-g^{\mu^*}_\infty>l^*$, so
there exists $T^*>0$ such that $h^{\mu^*}(T^*)-g^{\mu^*}(T^*)>l^*$. Applying
the continuous dependence of $(I_b , I_m ; g, h)$ on $\mu$, we can find
$\epsilon>0$ sufficiently small such that the solution of  (\ref{a3}) denoted by
$(I_{b}^{\mu}, I_{m}^{\mu}; g^\mu, h^\mu)$ satisfies
$$h^\mu (T^*)-g^\mu (T^*)>l^*$$
for all $\mu^* -\epsilon\leq\mu\leq\mu^*+\epsilon$. It follows that, for $\mu^* -\epsilon\leq\mu\leq\mu^*+\epsilon$,
$$h_\infty ^\mu-g_\infty ^\mu:=\lim_{t\to \infty} \ (h^\mu (t)-g^\mu (t))>l^*,$$
which implies that spreading
happens for $\mu\in [\mu^* -\epsilon, \mu^*+\epsilon]$ and therefore contradicts the definition of $\mu^*$. Hence $\mu^*\in\sum$.  The proof is completed.
 \epf

\bigskip
Next, we want to know what is the natural tendency of the virus when spreading happens, and therefore study the asymptotic behavior of the solution to problem \eqref{a3}.
\begin{thm} Assume that $R_0>1$. If spreading occurs, then the solution to the
 free boundary problem \eqref{a3} satisfies $\lim_{t\to +\infty} \ (I_b (x,t),I_m (x,t))=(I_b ^*, I_m ^*)$
uniformly in any bounded subset of $(-\infty, \infty)$, where $(I_b ^*, I_m ^*)$ is the unique positive equilibrium of system
\eqref{aode1}.
\end{thm}
\bpf  For clarity, we divide the proof into three steps.

(1) The superior limit of the solution

According to the comparison principle, we have $(I_b (x,t),I_m (x,t))\leq (\overline I_b (t), \overline I_m (t))$
for $g(t)\leq x\leq h(t)$, $t\geq 0$, where
$(\overline I_b (t), \overline I_m (t))$ is the solution of the problem
    \begin{eqnarray}
\label{ode1}
\left\{
\begin{array}{lll}
&\overline I_b '(t)= -\gamma_{b}\overline I_b (t)+ \alpha_{b} \beta_{b}\frac{(N_b -\overline I_b (t))}{N_b} \overline I_m (t),&t>0, \\
&\overline I_m '(t)= -d_{m}\overline I_m (t)+ \alpha_{m} \beta_{b} \frac{(A_m -\overline I_m (t))}{N_b} \overline I_b (t),&t>0,\\
&\overline I_b (0)=||I_{b,0}||_{L^\infty[-h_0, h_0]},\ \overline I_m (0)=||I_{m,0}||_{L^\infty[-h_0, h_0]}.&
\end{array} \right.
\end{eqnarray}
Since $R_0>1$, the unique positive equilibrium $(I_b ^*, I_m ^*)$ is globally asymptotically stable for the ODE system (\ref{ode1}) and $\lim_{t\to\infty}(\overline
I_b (t), \overline I_m (t))= (I_b ^*, I_m ^*)$; hence we obtain
\begin{equation}\label{123}
\limsup_{t\to +\infty} \ (I_b (x,t), I_m (x,t))\leq (I_b ^*, I_m ^*)
\end{equation}
uniformly for $x\in (-\infty, \infty)$.

(2) The lower bound of the solution for a large time

It is clear that
$$\lim_{l\to \infty} \frac{\alpha_{m} \alpha_{b} \beta_{b}^2 A_m}{N_{b} d_{m}(\gamma_{b}+D_1 (\frac \pi {2l})^2)}=R_0 >1,$$
we then can select some $L_0>0$ such that $\frac {\alpha_{m} \alpha_{b} \beta_{b}^2 A_m}{N_{b} d_{m}(\gamma_{b}+D_1 (\frac \pi {2L_0})^2)}>1.$
This implies that the principal eigenvalue $\lambda_0^*$ of
\begin{eqnarray}
\left\{
\begin{array}{lll}
-D_1 \psi_{xx}=-\gamma_{b} \psi+\frac {\alpha_{m} \alpha_{b} \beta_{b}^2 A_m}{N_{b} d_{m}}\psi+\lambda_0^* \psi,\; &
x\in (-L_0, L_0),  \\
\psi(x)=0, &x=\pm L_0
\end{array} \right.
\label{B1k}
\end{eqnarray}
satisfies
$$\lambda_0^*=\gamma_{b}+D_1 (\frac \pi {2L_0})^2-\frac {\alpha_{m} \alpha_{b} \beta_{b}^2 A_m}{N_{b} d_m}<0.$$
When the spreading happens, $h_\infty-g_\infty=\infty$, and then $h_\infty=-g_\infty=\infty$ from Lemma 3.1. Therefore,
 for any $L\geq L_0$, there exists $t_L>0$ such that $g(t)\leq -L$ and $h(t)\geq L$ for $t\geq t_L$.
Taking $\underline u=\delta \psi$ and $\underline v=\frac{\alpha_{m} \beta_{b} A_m}{N_b d_m} \underline u$, we can choose $\delta$ sufficiently small such that $(\underline u, \underline v)$
satisfies
 \begin{eqnarray*}
\left\{
\begin{array}{lll}
\underline u_{t}\leq D_1 \underline u_{xx}-\gamma_{b}\underline u+\alpha_{b} \beta_{b}\frac{(N_b -\underline u)}{N_b} \underline v,\; &  -L_0<x<L_0,\ t>t_{L_0},  \\
\underline v_{t}\leq-d_{m}\underline v+ \alpha_{m} \beta_{b} \frac{(A_m -\underline v)}{N_b} \underline u,\; & -L_0<x<L_0,\ t>t_{L_0},  \\
\underline u(x,t)=\underline v(x,t)=0,\ \quad & x=\pm L_0,\ t>t_{L_0},\\
\underline u(x,t_{L_0})\leq I_b (x,t_{L_0}),\ \underline v(x,t_{L_0})\leq I_m (x, t_{L_0}),
 & -L_0\leq x\leq L_0,
\end{array} \right.
\end{eqnarray*}
 this implies that $(\underline u, \underline v)$ is a lower solution of the solution $(I_b ,I_m)$ in $[-L_0, L_0]\times [t_{L_0}, \infty)$.
We then have $(I_b , I_m)\geq (\delta \psi, \frac{\alpha_{m} \beta_{b} A_m}{N_b d_m} \delta \psi)$ in $[-L_0, L_0]\times [t_{L_0}, \infty)$, which means that the solution can not tends to zero.

(3) The inferior limit

We first extend $\psi(x)$ to $\psi_{L_0}(x)$ by letting $\psi_{L_0}(x):=\psi(x)$ for $-L_0\leq x\leq L_0$ and $\psi_{L_0}(x):=0$ for $x<-L_0$ or $x>L_0$.
For $L\geq L_0$, $(I_b, I_m)$ satisfies
\begin{eqnarray}
\left\{
\begin{array}{lll}
\frac{\partial I_{b}}{\partial t}=D_1 \frac{\partial I_{b}}{\partial x^2}-\gamma_{b}I_b +\alpha_{b} \beta_{b} \frac{(N_b -I_b)}{N_b} I_m,\; &  g(t)<x<h(t),\ t>t_L,  \\
\frac{\partial I_{m}}{\partial t}=-d_{m} I_m + \alpha_{m} \beta_{b} \frac{(A_m-I_m)}{N_b} I_b,\; & g(t)<x<h(t),\ t>t_L,  \\
I_b (x,t)=I_m (x,t)=0, \quad & x=g(t)\, \textrm{or}\, x=h(t),\ t>t_L,\\
I_b (x,t_L)\geq \delta \psi_{L_0},\ I_m (x,t_L)\geq \frac{\alpha_{m} \beta_{b} A_m}{N_b d_m} \delta \psi_{L_0},
 & -L\leq x\leq L;
\end{array} \right.
\label{fs1}
\end{eqnarray}
therefore we have $(I_b , I_m)\geq (\underline I_b , \underline I_m)$ in $[-L, L]\times [t_L, \infty)$,
where $(\underline I_b , \underline I_m)$ satisfies
\begin{eqnarray}
\left\{
\begin{array}{lll}
\frac{\partial \underline I_{b}}{\partial t}\ = D_1 \frac{\partial \underline I_{b}}{\partial x^2}-\gamma_{b}\underline I_b +\alpha_{b} \beta_{b}\frac{(N_b - \underline I_b)}{N_b} \underline I_m ,\; &  -L<x<L,\ t>t_L,  \\
\frac{\partial\underline I_{m}}{\partial t} = -d_{m}\underline I_m +\alpha_{m} \beta_{b} \frac{(A_m -\underline I_m)}{N_b} \underline I_b,\; & -L<x<L,\ t>t_L,  \\
\underline I_b (x,t)=\underline I_m (x,t)=0,\quad & x=\pm L,\ t>t_L,\\
\underline I_b (x,t_L)=\delta \psi_{L_0},\ \underline I_m (x,t_L)= \frac{\alpha_{m} \beta_{b} A_m}{N_b d_m}\delta \psi_{L_0},
 & -L\leq x\leq L.
\end{array} \right.
\label{fs11}
\end{eqnarray}
It is easy to see that the model (\ref{fs11}) is quasimonotone increasing,
according to the upper and lower solution method
 and the theory of monotone dynamical systems (\cite{HS}, Corollary 3.6), we deduce that
$\lim_{t\to +\infty} \ (\underline I_b (x,t), \underline I_m (x,t))\geq (\underline I_{b,L}(x), \underline I_{m,L}(x))$ uniformly in
$[-L, L]$, where $(\underline I_{b,L}, \underline I_{m,L})$ satisfies
\begin{eqnarray} \label{fs12}\left\{
\begin{array}{lll}
-D_1 \underline I''_{b,L}=-\gamma_{b}\underline I_{b,L}+\alpha_{b} \beta_{b}\frac{(N_b -\underline I_{b,L})}{N_b} \underline I_{m,L},\; &  -L<x<L,  \\
-d_{m}\underline I_{m,L}+\alpha_{m} \beta_{b} \frac{(A_m -\underline I_{m,L})}{N_b} \underline I_{b,L}=0,\; & -L<x<L,  \\
\underline I_{b,L}(x)=\underline I_{m,L}=0, &x=\pm L
\end{array} \right.
\end{eqnarray}
and it is  the minimal upper solution over $(\delta \psi_{L_0},\ \frac{\alpha_{m} \beta_{b} A_m} {N_b d_m}\delta \psi_{L_0})$.

It follows from the comparison principle that the solution is increasing with $L$, that is, if $0<L_1<L_2$, then $(\underline I_{b, L_1}(x),\underline I_{m, L_1}(x)) \leq
(\underline I_{b, L_2}(x), \underline I_{m, L_2}(x))$
 in $[-L_1, L_1]$. Letting $L\to \infty$ and applying a classical elliptic regularity theory and a diagonal procedure yield that  $(\underline I_{b,L}(x), \underline I_{m,L}(x))$ converges
uniformly on any compact subset of $(-\infty, \infty)$ to $(\underline I_{b,\infty}, \underline I_{m,\infty})$, where $(\underline I_{b,\infty}$, $\underline I_{m,\infty})$ is continuous on $(-\infty, \infty)$ and satisfies
\begin{eqnarray*} \left\{
\begin{array}{lll}
-D_1 \underline I''_{b,\infty}=-\gamma_{b}\underline I_{b,\infty}+\alpha_{b} \beta_{b} \frac{(N_b -\underline I_{b,\infty})}{N_b} \underline I_{m,\infty},\; &  -\infty<x<\infty,  \\
-d_{m}\underline I_{m,\infty}+\alpha_{m} \beta_{b} \frac{(A_m -\underline I_{m,\infty})}{N_b} \underline I_{b,\infty}=0,\; & -\infty<x<\infty,  \\
\underline I_{b,\infty}(x)\geq \delta \psi_{L_0},\ \underline I_{m,\infty}(x)\geq  \frac{\alpha_{m} \beta_{b} A_m}{N_b d_m} \delta \psi_{L_0}, &-\infty<x<\infty.
\end{array} \right.
\end{eqnarray*}

Now, we claim that $\underline I_{b,\infty} (x)\equiv I_b ^*$ and $\underline I_{m,\infty} (x)\equiv I_m ^*$.
In fact, the second equation shows that
$$\underline I_{m,\infty}=\frac{\alpha_{m}\beta_{b}A_{m}\underline I_{b,\infty}}{N_b d_m+\alpha_{m}\beta_{b}\underline I_{b,\infty}}.$$
which leads the first equation to become
$$-D_1 \underline I''_{b,\infty}=-\gamma_{b}\underline I_{b,\infty}+(\alpha_{b} \beta_{b} \frac{(N_b -\underline I_{b,\infty})}{N_b}) \frac{\alpha_{m}\beta_{b}A_{m}\underline I_{b,\infty}}{N_b d_m+\alpha_{m}\beta_{b}\underline I_{b,\infty}}.$$
Considering the problem
$$-D_1 u''=-\gamma_{b}u+(\alpha_{b} \beta_{b} \frac{(N_b -u)}{N_b}) \frac{\alpha_{m}\beta_{b}A_{m}u}{N_b d_m+\alpha_{m}\beta_{b}u}:=f(u)u.$$
on can easily see that $f(u)$ is decreasing, so the positive solution is unique and $\underline I_{b,\infty} (x)\equiv I_b ^*$, therefore, $\underline I_{m,\infty} (x)\equiv I_m ^*$.

Based on the above fact, for any given $[-N, N]$ with $N\geq L_0$, we have that $(\underline I_{b,L}(x), \underline I_{m,L}(x))\to (I_b ^*, I_m ^*)$ uniformly in $[-N, N]$ as $L\to \infty$, and for any $\varepsilon >0$, there exists $L^*>L_0$ such that  $(\underline I_{b, L^*}(x), \underline I_{m, L^*}(x))\geq  (I_b ^*-\varepsilon, I_m ^*-\varepsilon)$ in $[-N, N]$. As above, there is $t_{L^*}$ such that $[g(t), h(t)]\supseteq [-L^*, L^*]$ for $t\geq t_{L^*}$.
So,
$$(I_b (x,t), I_m (x,t))\geq (\underline I_b (x,t), \underline I_m (x,t))\ \textrm{in}\ [-L^*, L^*]\times [t_{L^*}, \infty),$$ and
$$\lim_{t\to +\infty} \ (\underline I_b (x,t), \underline I_m (x,t))\geq (\underline I_{b, L^*}(x), \underline I_{m, L^*}(x))\ \textrm{in}\ [-L^*, L^*],$$
which together with the fact that $(\underline I_{b,L^*}(x), \underline I_{m,L^*}(x))\geq (I_b ^*-\varepsilon, I_m ^*-\varepsilon)$ in $[-N, N]$ gives
 $$\liminf_{t\to +\infty} \ (I_b (x, t), I_m (x,t))\geq (I_b ^*-\varepsilon, I_m ^*-\varepsilon)\ \textrm{in}\ [-N, N].$$
Since $\varepsilon>0$ is arbitrary,  we have $\liminf_{t\to +\infty} \ I_b (x, t)\geq I_b ^*$ and $\liminf_{t\to +\infty} \ I_m (x,t)\geq I_m ^*$ uniformly in $[-N,N]$, which together with (\ref{123}) concludes that $\lim_{t\to +\infty} \ I_b (x,t)=I_b ^*$
and  $\lim_{t\to +\infty} \ I_m (x,t)=I_m ^*$ uniformly in any bounded subset of $(-\infty, \infty)$.

\epf

Naturally, we can obtain the following spreading-vanishing dichotomy theorem, after combining Remark 4.1, Theorems 4.3 and 4.5.
\begin{thm} Assume that $R_0>1$.
Let $(I_b (x, t), I_m (x, t); g(t), h(t))$ be the solution of  free boundary problem \eqref{a3}.
Thus, the following dichotomy holds:

Either
\begin{itemize}
\item[$(i)$] {\rm Spreading:} $h_\infty-g_\infty =+\infty$ and $\lim_{t\to +\infty} \ (I_b (x, t), I_m (x, t))=(I_b ^*, I_m ^*)$
uniformly in any bounded subset of $(-\infty, \infty)$; \end{itemize}

or
\begin{itemize}
\item[$(ii)$] {\rm Vanishing:} $h_\infty -g_\infty \leq l^*$ with $\frac {\alpha_{m} \alpha_{b} \beta_{b}^2 A_m}{N_{b} d_m (
\gamma_{b}+D_1 (\frac \pi {l^*})^2)}=1$ and $\lim_{t\to +\infty} \ (||I_b (\cdot, t)||_{C([g(t),h(t)])}+||I_m (
\cdot, t)||_{C([g(t),h(t)])})=0$.
\end{itemize}
\end{thm}

\section{Numerical simulation and discussion }

In this section, first we give some simulations to illustrate our analytical results.
We set the following constants as given in \cite{LRD}, namely, $ A_m/N_b=20, \alpha_b=0.88, \alpha_m=0.16$  and $d_m=0.029.$

 From  Fig. 1, one can see that the virus in a scenario of vanishing for some small $\mu(:=0.1$) with $R_0^F (0)<1 $; in this case the infected population of birds will tends to zero gradually and the free boundaries expand slowly.

\begin{figure}[ht]
\centering
\subfigure[]{{
\includegraphics[width=0.4\textwidth]{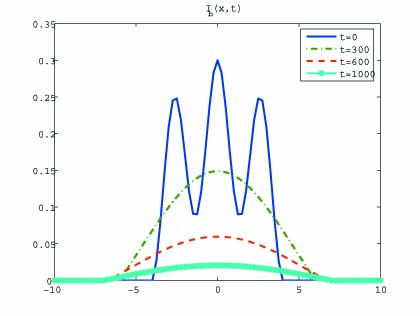}}}
\hspace{0.1cm}
\subfigure[]{{
\includegraphics[height=2in, width=0.4\textwidth]{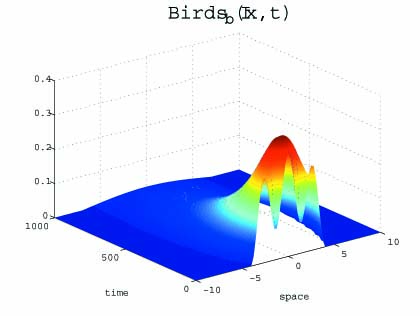}}}
\renewcommand{\figurename}{Fig.}
\caption{\scriptsize {Virus vanishing  ($R_0^F (0) < 1<R_0=1.31$); $D_1 = 4$; $\beta_b=0.09$; $\gamma_{b} = 0.6$
; $\mu=0.1$; $h_0=4.$}}
\end{figure}

Fig. 2 illustrates the spreading of the virus for $R_0^F (0)> 1 $; in this case,
 it is easy to see that the infected population of birds will not decay to zero and the free boundaries are expanding.
\begin{figure}[ht]
\centering
\subfigure[]{{
\includegraphics[width=0.4\textwidth]{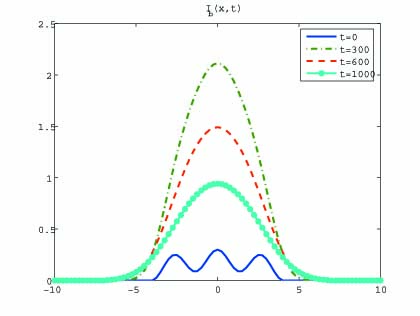}}}
\hspace{0.1cm}
\subfigure[]{{
\includegraphics[height=2in, width=0.4\textwidth]{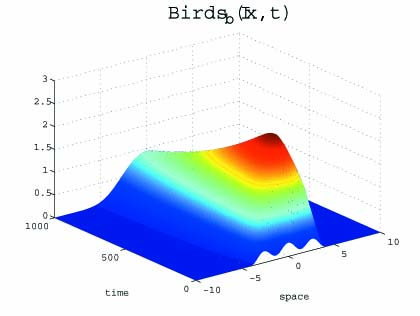}}}
\renewcommand{\figurename}{Fig.}
\caption{\scriptsize {Virus spreading  ($R_0^F (0)> 1 $);
$D_1 = 0.001$; $\beta_b=0.3$; $\gamma_{b} = 0.001$; $\mu=0.001$;
$h_0=4.$}}
\end{figure}

Comparing above two figures, we can see that when the expanding capacity $\mu$ is small and $R_0^F (0)< 1$, vanishing occurs, namely, the virus will be controlled (Fig. 1). On the other hand for any $\mu>0$ with  $R_0^F (0)> 1$, spreading happens; in this case the virus will spread gradually and the whole bounded area will infected by the virus in long run (Fig. 2).

 In attempt to describe the gradual spreading process and changing of the domain, free boundary problems, especially well-known Stefan conditions, have been used in mathematics and related disciplines.
 In this paper, we have examined the dynamic behavior of the populations $I_b (x,t)$ and $I_m (x,t)$ with double expanding fronts $x=g(t)$ and $x=h(t)$ modeled by system (\ref{a3}), which contains a coupled equations to describe the diffusion of birds by a PDE and the movement of mosquitoes by a ODE.

  We have presented the sufficient conditions for the WNv to be spreading or vanishing. Here, the vanishing means that the infected environment is limited and the virus disappears gradually (Fig. 1), while the spreading implies that the infected habitat is expanding to the whole environment and the virus always exists (Fig. 2).

  For the spatially-independent model (\ref{aode1}), it is shown that the virus vanishes eventually for any initial values if $R_0 \leq 1$  or remain epidemic if $R_0 >1$. For the diffusive model (\ref{Aa1}), we also introduce a threshold parameter $ R_0 ^D (:=\frac{\alpha_{m} \alpha_{b} \beta_{b}^2 A_m}{N_b d_m (\gamma_{b}+D_1 \lambda_1)})$,  such that for $0<R_0 ^D <1$, the epidemic eventually tends to vanishing, while for $R_0 ^D >1$  a spatially inhomogeneous stationary epidemic state appears and is globally asymptotically stable.
However, in the model \eqref{a3} with free boundary, the infected interval is changing with the time $t$, therefore, the basic reproduction number is not a constant and should change with $t$. So we here define it as the risk index $R_0 ^F (t) (:= \frac{\alpha_{m} \alpha_{b} \beta_{b} ^2 A_m}{N_b  d_m(\gamma_{b}+D_1 (\frac{\pi}{h(t)-g(t)})^2)})$, which depends on the habitat $(g(t),h(t))$, diffusion coefficient of birds $D_1$ and other parameters in (\ref{a3}).
 Our results show that if $R_0 \leq 1$ the virus always vanishes (Theorem 3.4), but if $ R_0 ^F (t_0) \geq 1$ for some $t_0 \geq 0 $, the virus is spreading (Theorem 4.1 and Remark 4.1). For the case $R_0 ^F (0)<1<R_0$, the spreading or vanishing of the virus depends on the initial number of infected birds, the expanding capacity $\mu$, the length of initial habitat, the diffusion rate of birds and other factors (Theorem 3.6 and Theorem 4.3).  Furthermore, the spreading-vanishing dichotomy is given and a sharp threshold related to the expanding capacity is also presented to distinguish the spreading and vanishing of WNv.

   Biologically, our model (\ref{a3}) is more realistic than the models (\ref{Aa1}) and (\ref{aode1}), because it gives a way to understand the diffusion process of infected birds and the movement of infected mosquitoes. Our theoretical results not only help us to understand which factors influence the spreading or vanishing of WNv, but also have useful implications for the control and elimination of WNv.


\end{document}